\documentclass{gtpart}     
\agtart
\usepackage{pinlabel}  
\usepackage{graphicx}  
\usepackage[all]{xy}
\title[Equivariant Towers and Non-positively Curved Groups]{Lifting Group Actions, Equivariant Towers and Subgroups of Non-positively Curved Groups}

\author[R.G.~Hanlon]{Richard Gaelan Hanlon}
\givenname{Richard Gaelan}
\surname{Hanlon}
\email{rgh262@mun.ca}

\author[E.~Mart\'inez-Pedroza]{Eduardo Mart\'inez-Pedroza}
\givenname{Eduardo}
\surname{Mart\'inez-Pedroza}
\address{Memorial University\\ St. John's, Newfoundland and Labrador, Canada}
\email{emartinezped@mun.ca}
\urladdr{www.math.mun.ca/~emartinezped/}


\arxivreference{math.GR/1307.2640}
\arxivpassword{sk7bs}

\volumenumber{}
\issuenumber{}
\publicationyear{}
\papernumber{}
\startpage{}
\endpage{}
\doi{}
\MR{}
\Zbl{}
\received{}
\revised{}
\accepted{}
\published{}
\publishedonline{}
\proposed{}
\seconded{}
\corresponding{}
\editor{}
\version{}

\newtheorem{theorem}{Theorem}[section]
\newtheorem{lemma}[theorem]{Lemma}
\newtheorem{proposition}[theorem]{Proposition}
\newtheorem{corollary}[theorem]{Corollary}
\newtheorem{remark}[theorem]{Remark}

\theoremstyle{definition}
\newtheorem{definition}[theorem]{Definition}

\newcommand {\N}{\mathbb{N}}

\newcommand {\Gtwo}{\ensuremath{\textit{G-2}}}
\newcommand {\Htwo}{\ensuremath{\textit{H-2}}}

\newcommand{\nclose}[1]{\ensuremath{\langle\!\langle#1\rangle\!\rangle}}

\DeclareMathOperator{\area}{\mathsf{Area}}
\DeclareMathOperator{\diam}{\mathsf{diam}}

\DeclareMathOperator{\kernel}{\mathsf{Kernel}}

\DeclareMathOperator{\full}{\mathsf{Span}}

\DeclareMathOperator{\cat0}{{CAT}(0)}

\DeclareMathOperator{\catK}{{CAT}(\kappa)}

\begin{document}

\begin{abstract}    If $\mathcal C$ is a class of complexes closed under taking full subcomplexes and covers and $\mathcal G$ is the class  of groups admitting proper and cocompact actions on one-connected complexes in $\mathcal C$, then $\mathcal G$ is closed under taking finitely presented subgroups. As a consequence the following classes of groups are closed under taking finitely presented subgroups: groups acting geometrically on regular $\cat0$ simplicial complexes of dimension $3$, $k$--systolic groups for $k\geq 6$, and groups acting geometrically on $2$--dimensional negatively curved complexes.  We also show that there is a finite non-positively curved cubical $3$--complex which is not homotopy equivalent to a finite non-positively curved regular simplicial $3$--complex.  We include applications to relatively hyperbolic groups and diagramatically reducible groups. The main result is obtained by developing a notion of equivariant towers which is of independent interest.
\end{abstract}

\maketitle

\section{Introduction}

We show that some classes of non-positively curved groups are closed under taking finitely presented subgroups.  We assume all spaces are combinatorial complexes and all maps are combinatorial, see Section~\ref{def:comb-cplx} for a definition.  A complex is \emph{one-connected} if it is connected and simply-connected.  A subcomplex $K$ of $X$ is \emph{full} if for any cell $\sigma \subset X$,  $\partial \sigma \subset K$ implies $\sigma \subset K$.  An action of a group $G$ on a  space $X$ is \emph{proper} if for all compact subsets $K$ of $X$ there are finitely many group elements $g$ such that $K \cap g(K) \neq \emptyset$. The action is \emph{cocompact} if there is a compact subset $K$ of $X$ such that the collection $\{gK\colon g\in G\}$ covers $X$.  Our main result is the following.

\begin{theorem}\label{thm:principle-full}
Let $\mathcal C$ be a category of complexes closed under taking full subcomplexes and topological covers. Let $\mathcal G$ be the category of groups acting properly and cocompactly, by combinatorial automorphisms, on one-connected complexes in $\mathcal C$. Then $\mathcal G$ is closed under taking finitely presented subgroups.
\end{theorem}

Some words about the literature.  A theorem of Steve Gersten states that finitely presented subgroups of hyperbolic groups of cohomological dimension $\leq 2$ are hyperbolic~\cite{Ge96}. This is a result only for dimension $2$,  since Noel Brady exhibited a $3$--dimensional compact locally $\cat0$ cube complex with hyperbolic fundamental group  containing a finitely presented subgroup which is not hyperbolic~\cite{Br99}.  Another result on the positive side is that finitely presented subgroups of torsion free systolic groups are systolic, this was proved by Daniel Wise~\cite{Wi03sys}. The proof of Theorem~\ref{thm:principle-full} builds on ideas of Gersten, Howie, and Wise in~\cite{Ge96, Ho81, Wi03sys}, and is self contained in the paper. 

\subsection{Sample Applications}

A \emph{regular simplicial complex} is a piecewise Euclidean simplicial complex where each $1$--cell has unit length.  A result of Rena Levitt shows that the category of regular locally $\cat0$ simplicial complexes of dimension $3$ is closed under taking full subcomplexes~\cite{Le13}. It is immediate that this category is closed under taking covers.  Recall that a group action on a metric space is said to be \emph{geometric} if is proper, cocompact and by isometries. 

\begin{corollary}\label{cor:3simplicial}
If $G$ acts geometrically on a regular $\cat0$ simplicial complex of dimension $3$, then any finitely presented subgroup acts geometrically on a regular $\cat0$ simplicial complex of dimension $3$.
\end{corollary}

It is not known whether Corollary~\ref{cor:3simplicial} holds for higher dimensions. However, the proof presented here does not generalize since  Levitt has exhibited regular locally $\cat0$  simplicial complexes of dimension $\geq 4$ containing full subcomplexes which are not locally $\cat0$.

In~\cite{Br99}, Noel Brady constructed a compact $3$--dimensional locally $\cat0$  cubical complex $X$ such that $\pi_1(X)$ contains a finitely presented subgroup which does not admit a finite classifying space.  Since every compact locally $\cat0$ space is a classifying space for its fundamental group, the previous corollary together with Brady's example implies the following statement.

\begin{corollary}
There is a $3$--dimensional finite locally $\cat0$ cubical  complex  which is not homotopy equivalent to a $3$--dimensional finite  regular locally $\cat0$ simplicial complex.
\end{corollary}

Levitt proved that non-positively regular $\cat0$ simplicial complexes of dimension $3$ are systolic complexes~\cite{Le13-2},  for this larger class our main theorem also applies. 

The notion of \emph{simplicial nonpositive curvature} for simplicial complexes was introduced by Januszkiewicz and \'Swi\c atkowski in~\cite{JS06} and independently by Haglund~\cite{Ha03} as a combinatorial analog of nonpositive curvature metric conditions.  A simplicial complex $L$ is  \emph{flag} if any set of vertices  which are pairwise connected by $1$--cells of $L$, spans a simplex in $L$. A simplicial complex $L$ is \emph{$k$--large}, $k\geq 6$, if $L$ is flag and there are no embedded cycles of length $<k$  which are full subcomplexes of $L$. A simplicial complex $X$ is \emph{locally $k$--large} if links of all simplices in $X$ are $k$--large. The fact that the category of  locally $k$--large simplicial complexes is closed under taking covers and full-subcomplexes immediately follows from the definitions. A group is \emph{$k$--systolic} if it acts properly and cocompactly by simplicial automorphisms on a one-connected locally $k$--large simplicial complex.  Wise proved that finitely presented subgroups of torsion free $k$--systolic groups are $k$--systolic using a tower argument~\cite{Wi03sys}. Corollary~\ref{cor:3systolic}  extends Wise's result to include groups with torsion. Previously to this work, the statement of Corollary~\ref{cor:3systolic} was also proved by Ga\v{s}per Zadnik using different methods~\cite{PrZa}.

\begin{corollary}\label{cor:3systolic}
For $k\geq 6$, if $G$ is a $k$--systolic group then any finitely presented subgroup of $G$ is $k$--systolic. 
\end{corollary}

In~\cite{Osaj13}, Osajda introduced the notion of complexes with $SD^*_2(k)$ links for $k\geq 6$. There is proved that the class of fundamental groups of   compact complexes with $SD^*_2(k)$ links is closed under taking finitely presented subgroups~\cite[Thm. 8.7]{Osaj13}. His proof is a tower argument for compact complexes, and in particular it is shown that  $SD^*_2(k)$ complexes are closed under taking covers and full subcomplexes. It follows that Theorem~\ref{thm:principle-full} applies to the class of $SD^*_2(k)$ providing an extension of Osadja's result.

A $2$--complex $X$ is \emph{negatively curved}  if it satisfies one of the following conditions:
\begin{enumerate}
\item  (metric condition) there is $\kappa<0$ such that $X$  admits the structure of a locally $\catK$ $M_\kappa$--complex, or
\item  (conformal condition)  there is an assignment of a non-negative real number, called  an \emph{angle}, to each corner of each $2$--cell  such that the sum of the angles on an $n$--gon is strictly less than $(n-2)\pi$ and links of $0$--cells satisfy \emph{Gromov's link condition}: every non-trivial circuit in the link is of angular measure at least $2\pi$. 
\end{enumerate}
A group acting geometrically on a one-connected negatively curved $2$--complex is word hyperbolic~\cite{BrHa99, GrHypGps}. Gersten proved that finitely presented subgroups of fundamental groups of finite negatively curved $2$--complexes are word hyperbolic~\cite[Thm 2.1]{Ge96}. One easily verifies that the category of negatively curved $2$--complexes is closed under taking subcomplexes and covers.  The following corollary extends Gersten's result.

\begin{corollary}[Subgroups of $2$--Diml Negatively Curved Groups]\label{thm:goal-hyp}
Let $Y$ be a one-connected negatively curved proper and cocompact \Htwo--complex. If  $G\leq H$ is finitely presented then $G$ is word hyperbolic.
\end{corollary}

We remark that Corollary~\ref{thm:goal-hyp} would follow from Gersten's result if $H$ is known to be virtually torsion free. However, it is an open question whether hyperbolic groups are virtually torsion free~\cite{GrHypGps}.  We also obtain an analogous result for relatively hyperbolic groups stated below.  

A $G$--complex $X$ is \emph{almost proper} if $G$ acts properly on the complement of the zero-skeleton of $X$. Observe that a proper action is almost proper. A group is called \emph{slender} if all its subgroups are finitely generated.  For definitions of relatively hyperbolic groups and fine complexes we refer the reader to Section~\ref{subsec:relhyp}.

\begin{theorem}[Theorem~\ref{thm:rel-hyp}]\emph{[Subgroups of 2--Diml Relatively Hyperbolic Groups]}\label{thm:rel-hyp-intro}\newline
Let $Y$ be a one-connected negatively curved, fine,  almost proper and cocompact \Htwo--complex such that  $H$--stabilizers of cells are slender.  If $G\leq H$ is finitely presented then $G$ is  hyperbolic relative to a finite collection of $G$--stabilizers of cells of $Y$.
\end{theorem}

The term of diagramatically reducible complex, defined below, was introduced by Gersten in connection with the study of equations over groups~\cite{Ge87}. The notion was first used, with different names, by  Chiswell, Collins and Huebschmann~\cite{CCH81} and Sieradski~\cite{Si83}.  Recall that an \emph{immersion} is a locally injective map, and a \emph{near-immersion} is a map which is locally injective except at $0$--cells of the domain. A $2$--complex $X$ is \emph{diagramatically  reducible} if there are no near-immersions $C\to X$, where  $C$ is a cell structure for the $2$--dimensional sphere.  

Since the composition of a near-immersion followed by an immersion is a near-immersion, the category of diagramatically reducible complexes is closed under taking covers and subcomplexes.   This category of complexes includes  locally $\cat0$ $2$--complexes, certain classes of small cancellation complexes,  conformal negatively curved $2$--complexes, spines of hyperbolic knots, and non-positively curved square complexes to name a few examples. Recall that a proper $G$--complex $X$ is a model for $\underline{E}G$ if for every finite subgroup $F\leq G$, the fixed point set $X^F$ is contractible. The following theorem is proved in Section~\ref{sec:applications}. 

\begin{theorem}[Theorem~\ref{thm:EG}]\emph{[Diagramatically Reduced Groups]}\label{thm:EG-intro}
Let $X$ be a diagramatically reducible one-connected proper \Htwo--complex. If $G\leq H$ is finitely presented then $G$ admits a diagramatically reducible  $2$--dimensional cocompact  model for $\underline{E}G$.
\end{theorem}

\subsection{Equivariant towers and the proof of Theorem~\ref{thm:principle-full}}  Briefly, Theorem~\ref{thm:principle-full} is proved using an extension of the tower method to equivariant maps. There are several applications of towers in combinatorial group theory as illustrated in~\cite{Ge96, Ho81, Ho84, HuWi01, Wi03sys, Wi04} and the results of this paper.  Towers is a geometric technique from $3$--manifold topology  introduced by  Papakyriakopolous~\cite{Pa57}, and later brought to combinatorial group theory by Howie~\cite{Ho84}.  A combinatorial map $X\to Y$ between connected CW-complexes is a \emph{tower} if it can be expressed as a  composition of inclusions and covering maps. A \emph{tower lifting} $f'$ of $f$ is a factorization $f=g \circ f'$ where $g$ is a tower. The lifting $f'$ is \emph{trivial} if $g$ is an isomorphism and the lifting is \emph{maximal} if the only tower lifting of $f'$ is the trivial one. It is well known that if $X$ is a finite complex, then any combinatorial map $X\to Y$ admits a maximal tower~\cite[Lem. 3.1]{Ho81}.  A tower is called an $\mathcal F$--tower if it is a composition of covering maps and  inclusions of full-subcomplexes.

By a \emph{locally finite complex} we mean a complex such that every  closed cell intersects finitely many closed cells.

\begin{theorem}[Theorem~\ref{thm:maximal-towers-full}]\emph{[Maximal Equivariant $\mathcal F$--Towers]}~\label{thm:maximal-towers-intro}
Let ${\bf f}\colon X \to Y$ be a $G$--map.  If $X$ is one-connected and $G$--cocompact and $Y$ is locally finite, then $\bf f$ has a maximal $\mathcal F$--tower lifting $f =g  \circ  f'$ where $g$ and $f'$ are $G$--maps.
\end{theorem}

An analogous result to Theorem~\ref{thm:maximal-towers-intro} where $Y$ is not required to be locally finite and instead $f'$ is only a maximal tower lifting (not a maximal $\mathcal F$--tower lifting) also holds, see Theorem~\ref{thm:maximal-towers}. This slightly different result is relevant to applications such as Theorem~\ref{thm:rel-hyp-intro}.  

In~\cite{Wi03sys}, there is a  result similar to Theorem~\ref{thm:maximal-towers-intro} stating that if $X\to Y$ is a map between finite simplicial complexes then there is a maximal \emph{expanded}-tower lifting. This is a different class of towers, and neither result subsume the other. The class of expanded towers works well in the setting of $2$--skeletons of systolic complexes which Wise used to prove Corollary~\ref{cor:3systolic} in the torsion-free case. 

A consequence of Theorem~\ref{thm:maximal-towers-intro} is the following.

\begin{theorem}[Theorem~\ref{thm:main-full}]\emph{[Existence of Immersed Cocompact $\mathcal F$--Cores]}~\label{thm:F-core} \label{thm:ImmCocCores} 
If $Y$ is a one-connected, proper and locally finite  $H$--complex and $G\leq H$ is finitely presented, then there is a one-connected cocompact $G$--complex $X$ and a  $G$--equivariant $\mathcal F$--tower $X\to Y$.
\end{theorem}

The proof of Theorem~\ref{thm:F-core} relies on the construction of a one-connected and cocompact  $G$--complex $X_0$ together with a  $G$--map $f\colon X_0\to Y$; this construction uses the hypothesis that $Y$ is one-connected. Given such a map, since $X_0$ is one-connected, there is an equivariant  maximal $\mathcal F$--tower lifting $f =g  \circ  f'$ with $g\colon X \to Y$. The maximality of $f'$ implies that $X$ is one-connected, and together with $Y$ being locally finite, that $X$ is $G$--cocompact.  

We also prove a version of Theorem~\ref{thm:F-core} without the assumption that $Y$ is locally finite, but with the weaker conclusion that the $G$--map $X\to Y$ is only a tower,  see Theorem~\ref{thm:main}.   The proof of our main result, Theorem~\ref{thm:principle-full}, follows immediately from  Theorem~\ref{thm:F-core}. 
\begin{proof}[Proof of Theorem~\ref{thm:principle-full}]
Let $Y$ be a one-connected complex in $\mathcal C$, let $H$ be a group acting properly and cocompactly on $Y$, and let $G\leq H$ be a finitely presented subgroup. Since $Y$ is a proper and cocompact $H$--complex, it is locally finite.
By Theorem~\ref{thm:ImmCocCores}, there is a one-connected  cocompact $G$--complex $X$ and an $G$--equivariant $\mathcal F$--tower $X\to Y$. By equivariance, $X$ is also a proper $G$--complex. Since $\mathcal C$ is closed under taking full-subcomplexes and covers, it follows that $X$ is in $\mathcal C$.  Therefore $G$ is in $\mathcal G$.
\end{proof}

\subsection{Outline of the Paper:} The rest of the paper is organized in four sections.  Section~\ref{sec:2} contains a result that provides sufficient conditions to lift a group action on a space to an intermediate cover.   Section~\ref{sec:towers} contains the definition of equivariant  towers and the proof of the existence of maximal equivariant $\mathcal F$--towers. Section~\ref{sec:cores} contains the proof of the existence of immersed cocompact $\mathcal F$--cores.  The last section contains the proofs of Theorems~\ref{thm:goal-hyp} and~\ref{thm:EG-intro}. 

\subsection*{Acknowledgments}
Thanks to Daniel Wise and Mark Sapir for both suggesting the applications to systolic complexes, and Daniel Wise for sharing his unpublished preprint~\cite{Wi03sys} which motivated part of the work of this paper. We also thank Damian Osajda for useful comments. We also thank Saya Mart\'inez-Aoki for assistance during this work.  Last but not least, we thank the anonymous referee for feedback, several suggestions and corrections, and pointing out a serious mistake in an earlier version of the paper.  We also acknowledge funding by the Natural Sciences and Engineering Research Council of Canada NSERC. 

\section{Lifting Group Actions}\label{sec:2}

In this section, all spaces are topological spaces which are path-connected, locally path-connected, and semilocally simply-connected. These are the  the standard hypotheses for the existence of universal covers.  All  maps between spaces are  continuous.  For standard results on covering space theory we refer the reader to Hatcher's textbook on Algebraic Topology~\cite{AT-Hatcher}.

\begin{definition} Let $X$ be a $G$--space and let $Y$ be an $H$--space. A  map $f\colon X\to Y$ is \emph{equivariant} with respect to a group homomorphism $f_\#\colon G\to H$ if $f(g.x)=f_\#(g) f(x)$ for every $g\in G$ and $x \in X$. 
 \end{definition}

 If $Y$ is an $H$--space then the universal cover $\widetilde Y$ is naturally an $\widetilde H$--space  and the covering map $\widetilde Y \to Y$ is equivariant with respect to a natural group homomorphism $\widetilde H \to H$ as the theorem below states. The proof  is patterned after an argument by Bridson and Haefliger  in the context of complexes of groups~\cite[Chap. III.C 1.15]{BrHa99}. 

\begin{theorem}[Lifting an Action to the Universal Cover] \label{thm:lifting}  Let $Y$ be an $H$--space, and  let $\rho: \widetilde{Y} \rightarrow Y$ be the universal covering map. Then there is a group $\widetilde H$ and action $\widetilde H \times \widetilde Y \to \widetilde Y$ with the following properties.
\begin{enumerate}
\item \label{p1}There is an exact sequence of groups 
\[ 1 \longrightarrow \pi_1 Y \longrightarrow \widetilde{H}  \longrightarrow  H \longrightarrow 1 .\]
\item \label{p3} The covering map $\widetilde Y \to Y$ is equivariant with respect to $\widetilde{H}  \to H$.
\item \label{p2}The restriction of $\widetilde H \times \widetilde Y \to \widetilde Y$ to $\pi_1 Y \times \widetilde Y \to \widetilde Y$ is the standard action by deck transformations of $\pi_1 Y$ on $\widetilde Y$.
\item \label{p4}  For each $\tilde y \in \tilde Y$ mapping to $y
\in Y$, the homorphism $\widetilde H\to H$ restricts  to an isomorphism $\widetilde H_{\tilde y} \to H_{y}$ between the $\widetilde H$--stabilizer of $\tilde y$ and the $H$--stabilizer of $y$. 
\item \label{p5} If $G\times \widetilde Y \to \widetilde Y$ is an action satisfying the four analogous properties above, then there is an isomorphism $\Phi\colon G\to \widetilde H$ such that $g.\tilde y=\Phi(g).\tilde y$ for all $\tilde y\in \widetilde Y$.
\item \label{p6} If $f\colon X \rightarrow Y$ is equivariant with respect to $f_\# \colon G \rightarrow H$, and $X$ is one-connected, and  $\tilde{f} \colon X \rightarrow \widetilde{Y}$ is a lifting of $f$, then $\tilde{f}$ is equivariant with respect to a group homomorphism $\tilde{f} \colon G \rightarrow \widetilde{H}$ which lifts  $f\colon G\to H$. 
\end{enumerate}
\end{theorem}

It is an immediate corollary that there exists liftings to intermediate covers which are $H$--regular as defined below.   

\begin{definition}[$H$--regular covers]\label{def:regular-G-cover}
Let $Y$ be an $H$--complex. A covering  map $\hat Y\to Y$ is \emph{$H$--regular} if the composition $\pi_1 \widehat Y \to \pi_1 Y \to \widetilde H$ is a normal subgroup of $\widetilde H$.
\end{definition}

\begin{remark}
Let $Y$ be an $H$--complex. By definition, the universal cover of an $H$--complex is $H$--regular. Furthermore  any cover of $Y$ associated to a characteristic subgroup of $\pi_1 (Y)$ is $H$--regular.
\end{remark}

\begin{corollary}[Lifting an Action to an Intermediate Cover]\label{cor:lifting}
Let $Y$ be an $H$--space and let $\widehat Y \to Y$ be an $H$--regular cover. Then the quotient group $\widehat H =\widetilde H/\pi_1\widehat Y$ acts on $\widehat Y$ and the map $\widehat Y \to Y$ is equivariant with respect to $\widehat H \longrightarrow H$.
Furthermore, stabilizers of points are preserved in the sense that  if $\hat y \in \widehat Y$ maps to $y \in Y$ then the restriction $\widehat H_{\hat y} \to H_y$ is an isomorphism.
\end{corollary}

\begin{proof}[Proof of of Theorem~\ref{thm:lifting}] 
  Let $y_0$ be a point of $Y$, and recall that $\widetilde Y$ can be identified with the set 
\[ \widetilde Y = \{[c]\ |\ \text{$c$ is a path in $Y$ starting at $y_0$} \}  ,\]
where $[c]$ denotes the homotopy class of $c$  with respect to homotopies fixing the endpoints $c (0)$ and $c (1)$. The covering map $\widetilde Y \to Y$ is interpreted as sending $[c]$ to $c (1)$.  The action of $\pi_1 (Y, y_0)$ on $\widetilde Y$ is given by \[ \pi_1 (Y, y_0) \times \widetilde Y  \longrightarrow \widetilde Y,\quad [\gamma] \times [c] \mapsto [\gamma\ast c].\]
For details of this standard construction of the universal cover and the action of the fundamental group, we refer the reader to~\cite{AT-Hatcher}.   Let  \[ \widetilde{H} = \{ (h,[c] ) \colon  h \in H,\text{ $c$ a path in $Y$ from $y_0$ to $h.y_0$} \}.  \]  
The group operation on $\widetilde H$ is given by 
\[ \widetilde H \times \widetilde H \longrightarrow \widetilde H, \quad (h,[c]) \times (h', [c'] ) \mapsto (hh', [c \ast h.c']),\]
where as usual $\ast$ denotes concatenation of paths. Observe that the operation is well defined since  for any pair of paths $f\in [c]$ and $f'\in [c']$, the terminal point of $f$ equals the initial point of $hf'$, and  $c\ast h.c'$ is homotopic relative to endpoints to $f \ast h.f'$, and the terminal point of  $c\ast h.c'$ is $hh'.y_0$. To show that this operation endows $\widetilde H$ with a group structure is routine and it is left to the interested reader. The action of $\widetilde H$ on $\widetilde Y$ is given by 
\[ \widetilde H \times \widetilde Y \longrightarrow \widetilde Y, \quad (h, [c] ).[c']  = [c \ast h.c'] ,\]
and one easily verifies that it is a well-defined action. Now we verify the six properties.

\begin{proof}[Properties~\eqref{p1},~\eqref{p3} and~\eqref{p2}] Observe that the natural projection 
\[ \widetilde H \longrightarrow H, \quad  (h, [c]) \mapsto h \]
is a surjective homomorphism with kernel 
\[  \left\{ (1, [c]) \  |\ \text{$c$ is a closed path with $c(0)=c(1)=y_0$} \right \} \cong \pi_1 (Y, y_0).\]
 By definition the action $\widetilde H \times \widetilde Y \to \widetilde Y$ extends the action of $\pi_1(Y, y_0)$ on $\widetilde Y$.  
 To verify that $\widetilde Y\to Y$ is equivariant with respect to $\widetilde H \to H$, let $(h', [c']) \in \widetilde{H}$ and $[c]  \in \widetilde Y$ and observe that the terminal point of the path $(h',[c']).[c]$ equals the point $h'.c(1)$.
\end{proof}

\begin{proof}[Property~\eqref{p4}]
This statement follows from properties~\eqref{p1},~\eqref{p3} and ~\eqref{p2} as follows. 
Let $\tilde y \in \tilde Y$ and let $y\in Y$ be the image of $\tilde y$ by the covering map. By equivariance, the homomorphism $\widetilde H \to H$ maps the stabilizer $\widetilde H_{\tilde y}$ into $H_{y}$. Let $h\in H_y$ and let $\widetilde h \in \tilde H$ be an element mapping to $h$. Then there is $g\in \pi_1 Y$ such that $g.\tilde y = \tilde h.\tilde y$. It follows that $g^{-1}\tilde h \in \widetilde H_{\tilde y}$ and $g^{-1}\tilde h$ maps to $h$.  This shows that $\widetilde H \to H$ is surjective. For injectivity, if $\tilde h_1, \tilde h_2 \in \tilde H_{\tilde y}$ map to $h\in H_y$, then $\tilde h_1^{-1}\tilde h_2 \in \pi_1 Y \cap \widetilde H_{\tilde y}$. Since $\pi_1 Y$ acts freely on $Y$, it follows that $\tilde h_1=\tilde h_2$.
\end{proof}

\begin{proof}[Property~\eqref{p5}]
 The isomorphism is a consequence of the short five lemma. Suppose that $G$ is another group acting on $\widetilde Y$ and satisfying properties one to three; the fourth property is not needed as it is a consequence of the other three.  By property (1) of $G$, there is a short exact sequence
\[1\longrightarrow K \longrightarrow G \overset{\varphi}{\longrightarrow } H \longrightarrow 1,\]
where $K\cong \pi_1 (Y)$. Denote by $\rho$ the covering map $\widetilde Y \to Y$, and let $\tilde y_0\in \widetilde Y$ be such that $\rho (\tilde y_0)=y_0$. For $g\in G$, let $c_g$ denote a path in $Y$ starting at $y_0$ obtained by composing $\rho$ with a path in $\widetilde Y$ from $\tilde y_0$ to $g.\tilde y_0$; here we use that $G$ acts on $\widetilde Y$. Since $\widetilde Y$ is simply-connected, the homotopy class $[c_g]$ depends only on $g$; by property (2) of $G$, the pair $(\varphi (g), [c_g])$ is an element of $\widetilde H$. It follows that  there is a well defined  group homomorphism 
\[ \Phi\colon G \longrightarrow \widetilde H, \quad g\mapsto ( \varphi (g), [c_g]).\]

Observe that $\Phi (K)$ is a subset of the kernel of $\widetilde H \to H$ and therefore $\Phi$ satisfies the following commutative diagram
\begin{equation}\nonumber 
\xymatrix{ 
   1 \ar[r] & K \ar[r] \ar[d]^\Phi & G  \ar[r]^\varphi \ar[d]^\Phi & H \ar[d]^{\text{id}} \ar[r] &1 \\
  1 \ar[r] & \pi_1(Y, y_0)  \ar[r] & \widetilde{H}   \ar[r] & H \ar[r] &1.
}\end{equation}
We claim that the $\Phi\colon K \to \pi_1 (Y, y_0)$ is an isomorphism.  By property (3) of $G$, the action of $K$ on $\widetilde Y$ is the action of deck transformations of $\widetilde Y \to Y$. Therefore for $[c] \in \pi_1 (Y, y_0)$, there is $k\in K$  acting on $\widetilde Y$ as $(1, [c])$ does. Recall that $[c]$ is identified with $(1, [c]) \in \widetilde H$.  By definition of $\Phi$, we have that $\Phi (k) = (\varphi (k), [c])=(1, [c])$ and hence $\Phi\colon K \to \pi_1 (Y, y_0)$ is surjective. For injectivity, observe that if $k\in K$ and $\Phi (k)$ is trivial, then $k \in G_{\tilde y_0}$; then, by property (3) of $G$, $k$ is trivial. By the short five lemma applied to the commutative diagram above, $\Phi\colon G \to \widetilde H$ is an isomorphism.

Now we verify that the actions of $\widetilde H$ and $G$ are identical up to composing with $\Phi$, i.e, for every $g\in G$ and $\tilde y \in \widetilde Y$,  $g.\tilde y = \Phi (g).\tilde y$. Suppose that $\tilde{y}= [c']$. Then $g.\tilde y = [c_g\ast \varphi (g).c']$ since  $c_g\ast \varphi (g).c'$ lifts to a path from $\tilde y_0$ to $g.\tilde y$. Hence $g.\tilde y = [c_g\ast \varphi (g).c'] = (\varphi (g), [c_g] ).[c']$.
\end{proof}

\begin{proof}[Property~\eqref{p6}]
To simplify notation, we denote by $f$ the map $X\to Y$ and the group homomorphism $G \to H$.  Let $x_0\in X$ be such that $f(x_0)=y_0$. Since $X$ is simply connected, using the description of $\widetilde Y$ as a set of equivalence classes of paths, for any $x\in X$ we have that $\tilde f (x) = [f \circ c]$  where $c$ is any path in $X$ from $x_0$ to $x$. We show that $\tilde{f}\colon X\to \widetilde Y$ is equivariant with respect to the group homomorphism 
\[\tilde f \colon  G \rightarrow \widetilde{H},\quad  g \mapsto (f(g),  [f\circ c]).\] 
where $c$ is a path in $X$ from $x_0$ to $g.x_0$. Observe that $\tilde f$ is well defined as a map since $X$ is simply connected and $f\colon X\to Y$ is an equivariant map. To show that $\tilde f$ is a homomorphism is routine. To verify equivariance,  first let $x\in X$ and $g\in G$. Let $c$ be a path from $x_0$ to $g.x$, let $c'$ be a path from $x_0$ to $g.x_0$, and let $c''$ be a path from $x_0$ to $x$.  Since $X$ is simply connected, we have that  $[c' \ast g.c''] = [c]$ and hence 
\begin{equation} \nonumber \begin{split} 
\tilde{f}(g).\tilde{f}(x) & = (f (g),[f \circ c'] ).[f \circ c'']  = [(f \circ c') \ast (f (g).(f \circ c''))]  \\ 
 & =[f(c' \ast g.c'')] = [f(c)]  = \tilde{f}(g.x). \qedhere \end{split} \end{equation}
\end{proof} This completes the proof of the theorem.
\end{proof}

\section{Maximal Equivariant Towers}\label{sec:towers}

For the rest of the paper, all spaces are combinatorial complexes and all maps are combinatorial.   All group actions on complexes are by combinatorial maps.
A $G$--complex $X$ is {\em proper} (respectively {\em cocompact, free}) if the $G$-action is proper (respectively, cocompact, free).  For a   cell $\sigma$ of $X$, the pointwise $G$--stabilizer of $\sigma$ is denoted by $G_\sigma$, and the $G$-orbit of $\sigma$ is denoted by $G(\sigma)$.

\begin{definition}[Combinatorial Complexes and Maps]\label{def:comb-cplx}\cite[Ch.I Appendix]{BrHa99}
A map $X\to Y$ between CW-complexes is combinatorial if its restriction to each open cell of $X$ is a homeomorphism onto an open cell of $Y$. A CW-complex $X$ is combinatorial provided that the attaching map of each open cell is combinatorial for a suitable subdivision. 
\end{definition}

\subsection{Definition of Equivariant Tower and Basic Lemmas}

\begin{definition}[Equivariant Map] 
Let $X$ be a $G$--complex and let $Y$ be an $H$--complex. An \emph{equivariant map ${\bf f}\colon X \to Y$} is a pair $(f, f_\#)$ where $f\colon X\to Y$ is a combinatorial map,  $f_\#\colon G\to H$ is a group homomorphism, and  $f$ is equivariant with respect to $f_\#$, that is,  $f(g.x)=f_\#(g) f(x)$ for every $g\in G$ and $x \in X$.  As usual, if  both $X$ and $Y$ are $G$-complexes and $f_\#$ is the identity map on $G$, then ${\bf f}$ is called a \emph{$G$--map}.  As a convention, we use bold letters to denote equivariant maps ${\bf f}$, and in this case, we use $f$ to denote the map between complexes, and $f_\#$ to denote the group homomorphism. The domain of an equivariant map is the space together with the group acting on it, and the same convention applies to the codomain. The composition of equivariant maps is defined in the natural way. Equality of equivariant maps ${\bf f}={\bf g}$ means their domains and codomains are equal, $f=g$, and $f_\#=g_\#$.
\end{definition}

\begin{definition}[Equivariant Isomorphism]
An equivariant map is an \emph{isomorphism} if the map at the level of spaces is a homeomorphism  and the map a the level of groups is a group isomorphism. 
\end{definition}

\begin{definition}[Equivariant Inclusions]
An equivariant map $\boldsymbol \imath = (\imath, \imath_\#)$ is called an \emph{equivariant inclusion} if  
$\imath$ and $\imath_\#$ are injective.  The equivariant inclusion $\boldsymbol \imath$ is called \emph{proper} if either $\imath$ or $\imath_\#$ is not surjective.
\end{definition}

\begin{definition}[Equivariant Cover]
Let $Y$ be an $H$--complex. A covering map $\widehat Y \to Y$ is called an \emph{equivariant cover} if $\widehat Y$ is an $H$--regular cover.  By Corollary~\ref{cor:lifting}, if $\rho$ is an equivariant cover then there is a well defined group homomorphism $\rho_\#\colon \widehat H \to H$ such that $(\rho, \rho_\#)$ is an equivariant map. When refering to an equivariant cover we will be implicitly referring to the associated equivariant map $(\rho, \rho_\#)$. 
\end{definition}

\begin{definition}[Towers and $\mathcal F$--Towers]\label{defn:equivariant-tower} Let $X$ be a $G$--complex and   let $Y$ be an $H$--complex. An  equivariant map ${\bf g}\colon X\to Y$ is an \emph{equivariant tower} if it can be expressed as an alternating composition of equivariant inclusions and equivariant covers. Specifically, $g$ is a composition
\begin{equation*} X=X_n \hookrightarrow \widehat X_{n} \to X_{n-1} \hookrightarrow \cdots \hookrightarrow \widehat X_2  \to X_1 \hookrightarrow \widehat X_1 \to X_0 = Y \end{equation*}
and $g_\#$ is a composition
\begin{equation*} G=G_n \hookrightarrow \widehat G_{n} \to G_{n-1} \hookrightarrow \cdots \hookrightarrow \widehat G_2  \to G_1 \hookrightarrow \widehat G_1 \to G_0 = H \end{equation*}
where  $X_i$ is a $G_i$-complex,  and
$\widehat X_{i} \to X_{i-1}$ is a $G_{i-1}$-regular cover of $X_{i-1}$ inducing the $\widehat G_i$-action on  $\widehat X_{i}$, and $X_{i}$ is a subcomplex of $\widehat X_i$ invariant under the subgroup $G_{i}\leq  \widehat G_i$, and both $G=G_n$, and $H=G_0$. In the case that each  $X_i$ is a full-subcomplex of $\widehat X_i$, the tower $\bf g$ is called an  \emph{$\mathcal F$--tower}.

The \emph{length} of the tower ${\bf g}$  is the smallest value of $n$ in an expression for ${\bf g}$ as above. In particular, an equivariant inclusion or cover have length at most one. By convention, the identity map has length zero.
\end{definition}

\begin{definition}[Equivariant Tower Lifting and $\mathcal F$--tower Lifting]
Let $\bf f$ be an equivariant map. An \emph{equivariant tower lifting} of ${\bf f}$ is an equivariant map $\bf f'$ such that there is an equivariant tower $\bf g$  such that ${\bf f}={\bf g}\circ {\bf f'}$. The lifting is \emph{trivial} if ${\bf g}$ is an equivariant isomorphism, and the lifting is \emph{maximal} if the only equivariant  tower lifting of ${\bf f'}$ is the trivial one.  The notions of 
\emph{equivariant $\mathcal F$--tower lifting}, \emph{trivial $\mathcal F$--tower lifting}, and \emph{maximal $\mathcal F$--tower lifting} are defined analogously.
\end{definition}

\begin{remark}[Composition of Towers]\label{rem:tower-comp}
Observe that if $\bf f$ and $\bf h$  are equivariant towers, and the codomain of $\bf f$ equals the domain of $\bf h$ (this means on the space and the group), then the composition ${\bf h} \circ {\bf f}$ is an equivariant tower. The same statement holds for  $\mathcal F$--towers.
 
Suppose ${\bf f}={\bf g}\circ {\bf f'}$ is a tower lifting of ${\bf f}$, and ${\bf f'}={\bf g'}\circ {\bf f''}$ is a tower lifting of ${\bf f'}$.  Since the composition ${\bf g} \circ {\bf g'}$ is a tower,  ${\bf f}=({\bf g} \circ {\bf g'}) \circ {\bf f''}$ is a tower lifting of ${\bf f}$. In particular, if ${\bf f''}$ is a maximal tower lifting of ${\bf f'}$ then ${\bf f''}$ is a maximal tower lifting of $\bf f$. The same statement holds for $\mathcal F$--towers liftings.
\end{remark}

\begin{definition}[$0$--surjective]
A map $X\to Y$ is \emph{$0$--surjective} if every $0$--cell of $Y$ is in the image of $X$.
\end{definition}

\begin{proposition}[Maximality $\Leftrightarrow$ Surjectivity]  \label{prop:maxitower-full} \label{prop:maxitower}
Let $X$ be a one-connected and $G$--cocompact complex,  let ${\bf f}\colon X\to Y$ be a $G$--map.  
\begin{enumerate}
\item\label{item:maxitower} An equivariant tower lifting ${\bf  f'}$ of $\bf f$ is maximal if and only if $f'$ is surjective and $\pi_1$--surjective, and $f'_\#$ is surjective.
\item\label{item:maxitower-full} An equivariant $\mathcal F$--tower lifting ${\bf  f'}$ of $\bf f$ is maximal if and only if $f'$ is $0$--surjective and $\pi_1$--surjective, and $f'_\#$ is surjective. 
\end{enumerate}
\end{proposition}
\begin{proof}
We sketch the proof for $\mathcal F$--towers and the proof of the first statement is left to the reader. 
The only if part is immediate since otherwise $\bf f'$ would factor through a non-trivial inclusion or the universal covering map. 
For the if part, let ${\bf f'} = {\bf h} \circ {\bf  f''}$ be  an $\mathcal F$--tower lifting of $\bf f'$. Suppose that ${\bf h}$ is an equivariant inclusion of a full-subcomplex; since $f'$ is $0$--surjective and $f'_\#$ is surjective we have that $h$ and $h_\#$ are surjective and hence $\bf h$ is an equivariant isomorphism. Suppose that ${\bf h}$ is an equivariant cover; since $X$ is one-connected and $f'\colon X \to Y'$ is $\pi_1$--surjective it follows that $Y'$ is one-connected and hence ${h}$ is the trivial cover and $h_\#$ is an isomorphism.  The general case follows by induction on the length of the tower.
\end{proof}

\begin{definition}[Preserving $0$--Stabilizers]
An equivariant map ${\bf f}$ from the $G$--complex $X$ to the $H$--complex $Y$ is said to \emph{preserve $0$--stabilizers} if for every $0$--cell $\sigma$ of $X$ the map $f_\#\colon G_\sigma \to H_{f(\sigma)}$ is a group isomorphism.
\end{definition}

\begin{definition}[$d({\bf f}), r({\bf f}), e(X)$]
For a $G$--complex $X$, let $v(G, X)$ denote the number of $G$-orbits of $0$--cells and let $e(G, X)$ denote the number of $G$-orbits of $1$--cells.  If the group is understood, we simply write $v(X)$ and $e(X)$.  Observe that  $v (G, X)<\infty$ and $e(G,X)<\infty$  if $X$ is $G$--cocompact.  If ${\bf f}\colon X \to Y$ is an equivariant map from the $G$--complex $X$ to the $f_\#(G)$-complex $Y$, we define $d({\bf f}) = v(G, X)$ and $r({\bf f})=v (f_\#(G), f(X))$. Observe that for a $G$--map ${\bf f}\colon X \to Y$ we have that $d({\bf f})\geq r({\bf f})$ and   $e(X)\geq e(f(X))$.
\end{definition}

\begin{lemma}\label{lem:injectiveOrbits-2} \label{lem:isomorphism}
Let ${\bf g}\colon X\to Y$ be an equivariant tower such that  $g_\#$ is an isomorphism,  preserves $0$--stabilizers, and $\infty>d({\bf g})=r({\bf g})$. The following statements hold.
\begin{enumerate}
\item \label{lem:t2r2} The map $\bf g$ is an equivariant inclusion and  when restricted to the $0$--skeletons is an isomorphism. 
\item  \label{lem:third} If $g$ is surjective, then $\bf g$ is an equivariant isomorphism.
\item \label{lem:ThanksToReferee} If  $\infty>e(X)=e(Y)$, then $g$ is an isomorphism when restricted to  $1$--skeletons and, in particular, $g$ is $\pi_1$--surjective.
\end{enumerate}
\end{lemma}
\begin{proof} Since $g_\#$ is an isomorphism, assume that  $X$ and $Y$ are  $G$-complexes and ${g_\#}\colon G\to G$ is the identity map.  Observe that the second statement is immediate if the first statement holds.

For the first statement, we show first that $\bf g$ induces a bijection between the $0$--skeletons. Since $\infty>d({\bf g})=r({\bf g})$,  
the pigeon-hole argument shows that ${\bf g}$ induces  a bijection between $G$-orbits of $0$--cells. Hence it is enough to show that for any $0$--cell $x$ in the domain of $\bf g$, the induced map between orbits $g\colon G(x) \to G(g(x)) $ is a bijection. By equivariance,  $g\colon G(x) \to G(g(x)) $ is surjective.  For injectivity, suppose there are $x, x'\in X$ and $h \in G$ such that $h.x=x'$ and $g(x)=g(x')$. By equivariance, $h$ is in the $G$--stabilizer of $g(x)$.  Since  $g_\#\colon G_{x} \to G_{g(x)}$ is the identity map, it follows that $h \in G_{x}$ and hence $x'=x$. We have proved that $g$ induces a bijection between the $0$--skeletons. Now we show that ${g}$ is an inclusion of complexes. Let $\sigma$ and $\sigma'$ be two $k$-cells of $X$ mapping to the same $k$-cell of $Y$.  Since $g$ is bijective on $0$--cells, the cells $\sigma$ and $\sigma'$ have a common $0$--cell in their closure. Since $g$ is a tower, it is a locally injective map. Therefore $\sigma$ and $\sigma'$ are the same $k$-cell.

The third statement is proved as follows. The first statement of the lemma implies that  $\bf g$ induces an equivariant inclusion between the $1$--skeletons, and an isomorphism between $0$--skeletons.  It remains to prove that the induced map between $1$--skeletons is surjective.  Since $\infty > e(X) = e(Y)$, the pigeon-hole argument shows that ${\bf g}$ induces  a bijection between $G$-orbits of $1$--cells.  For any $1$--cell $\sigma$, equivariance implies that the induced map between orbits $g\colon G(\sigma) \to G(g\sigma)$ is surjective. Therefore $g$ is surjective on  $1$--cells, and hence the induced  map between the $1$--skeletons is an isomorphism.  
\end{proof}

\subsection{Existence of Maximal $\mathcal F$--tower  Liftings}

\begin{theorem}[Maximal Equivariant $\mathcal F$--Towers]\label{thm:maximal-towers-full}
Let ${\bf f}\colon X \to Y$ be a $G$--map.  If $X$ is one-connected and $G$--cocompact and $Y$ is locally finite, then $\bf f$ has a maximal equivariant $\mathcal F$--tower lifting. 
\end{theorem}

Before the proof of the theorem we need a definition and a remark.

\begin{definition}[Span]
The \emph{span} of a subcomplex $K\subset X$,  denoted by $\full_X (K)$, is the smallest full subcomplex of $X$ containing $K$.
\end{definition}

\begin{remark}
If $Y$ is a $G$--complex and $K \subset Y$ is a $G$--subcomplex, then $\full_X (K)$ is a $G$--subcomplex of $Y$. If, in addition, $Y$ is locally finite and $K$ is $G$--cocompact then $\full_X (K)$ is $G$--cocompact and in particular $e(\full_X(K))$ is finite. 
\end{remark}

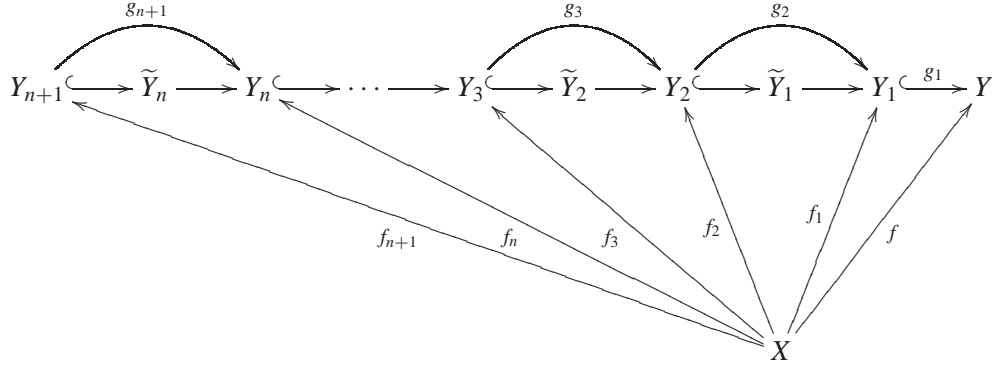
\begin{figure}
\[ \xymatrix@C=2.06em{   
 Y_{n+1} \ar@/^2pc/[rr]^{g_{n+1}} \ar@{^{(}->}[r] & \widetilde Y_n \ar[r]  &Y_n \ar@{^{(}->}[r] &  \cdots \ar[r] & Y_3  \ar@/^2pc/[rr]^{g_{3}} \ar@{^{(}->}[r] & \widetilde Y_2\ar[r] & Y_2\ar@{^{(}->}[r] \ar@/^2pc/[rr]^{g_{2}} &  \widetilde Y_1\ar[r] &  Y_1\ar@{^{(}->}[r]^{g_1}  & Y\\ \\ \\
  & &						 & 							 & 						   	   & 						    &  						      &  X \ar[uuurr]_{f} \ar[uuur]^{f_1} \ar[uuul]^{f_2} \ar[uuulll]^{f_3}   \ar[uuulllll]^{f_n}  \ar[uuulllllll]^{f_{n+1}}  \\
}
\]
\caption{The tower construction in the proof of Theorem~\ref{thm:maximal-towers-full}}\label{fig:tower-construction}
\end{figure}

\begin{proof}[Proof of Theorem~\ref{thm:maximal-towers-full}]
Let ${\bf f}_0$ denote the $G$--map ${\bf f}$. Let $Y_1$ be the $G$--subcomplex $\full ( f (X) )$ of $Y$,  let ${\bf f}_1\colon X \to  Y_1$ by the $G$--map defined as ${\bf f}_1(x)={\bf f}(x)$ for $x\in X$, and let ${\bf g}_1\colon Y_1 \to Y$ be the $G$--equivariant inclusion. Since $X$ is $G$--cocompact, $\infty>d({\bf f_1})-r({\bf f_1})$.  Since $Y$ is locally finite and $X$ is $G$--cocompact, $Y_1$ is locally finite and $\infty>e(Y_1)$.

For $n\geq 1$, suppose that for we have defined an $\mathcal F$--tower lifting ${\bf f}_n \colon X \to Y_n$ of ${\bf f}_{n-1}$ and an $\mathcal F$--tower ${\bf g}_n\colon Y_n \to Y_{n-1}$  such that ${\bf g}_n\circ {\bf f}_n = {\bf f}_{n-1}$ and $Y_n$ is locally finite and $\infty>d({\bf f}_n)-r({\bf f}_n)$ and $\infty>e(Y_n)$.   
Let ${\bf \rho}_n \colon \widetilde Y_n \to Y_n$ be the universal covering map and let  ${\bf \tilde f}_{n}\colon X\to \widetilde Y_n$ be the equivariant lifting of ${\bf f}_n$ to the universal cover $\widetilde Y_n$. By Theorem~\ref{thm:lifting}, ${\bf \tilde f}_{n}$ is naturally a $G$--map.   Let  $Y_{n+1}$ be the $G$--subcomplex $\full ( \tilde f_{n}(X) )$ of $\widetilde Y_n$ and let ${\bf f}_{n+1}\colon X \to Y_{n+1}$ the $G$--map defined by ${\bf f}_{n+1} (x) = {\bf \tilde f}_n(x)$ for each $x\in X$. Let   ${\bf g}_{n+1} \colon Y_{n+1} \to Y_n$ be the restriction of ${\bf \rho}_n$ to $Y_{n+1}$.  By construction, ${\bf f}_{n+1}\colon X \to Y_{n+1}$ is $0$--surjective. Since $\widetilde Y_n$ is locally finite and $X$ is $G$--cocompact, $Y_{n+1}$ is also locally finite and $\infty>e(Y_{n+1})$. The construction is illustrated in Figure~\ref{fig:tower-construction}. 

Consider the well-ordered set consisting of pairs of positive integers $\N \times \N$ with the dictionary order $\leq_{dic}$. For each ${\bf f}_n$ we assign an element $|{\bf f}_n|$ of  $\N \times \N$ defined as 
\[  |{\bf f}_n| = \Big(\ d({\bf f})-r({\bf f}_n)\ ,\ e(Y_n)\ \Big) .\]

\begin{lemma}[Decreasing complexity]
 The following statements hold.
\begin{enumerate}
\item If $Y_n$ is simply-connected then ${\bf f}_n\colon X\to Y_n$ is a maximal $\mathcal F$--tower lifting of ${\bf f}\colon X\to Y.$
\item If $|{\bf f}_{n+1}| =_{dic} |{\bf f}_{n}|$ then $Y_n$ is simply-connected.
\item For every $n$, $|{\bf f}_{n+1}| \leq_{dic} |{\bf f}_{n}|$. 
\end{enumerate}
\end{lemma}
\begin{proof}
If $Y_n$ is simply-connected then ${f}_n$ is  $\pi_1$--surjective. Since  ${f}_n$ is also $0$--surjective, Proposition~\ref{prop:maxitower-full} implies that ${\bf f}_n$ is a maximal $\mathcal F$--tower lifting of ${\bf f}$. This proves the first statement.

For the second statement, suppose that $|{\bf f}_{n+1}| =_{dic} |{\bf f}_{n}|$.  We will show that ${g}_{n+1}\colon Y_{n+1} \to Y_n$ is $\pi_1$--surjective using  Lemma~\ref{lem:isomorphism}\eqref{lem:ThanksToReferee}; then the proof concludes by observing that $g_{n+1}$ factors through the simply-connected space $\widetilde Y_n$ and hence $Y_n$ is simply-connected. It remains to verify that ${\bf g}_{n+1}$ satisfies the hypotheses of  Lemma~\ref{lem:isomorphism}\eqref{lem:ThanksToReferee}. By construction,  ${\bf g}_{n+1}\colon Y_{n+1} \to Y_n$ is a $G$--map and a tower; the assumption implies that  $\infty>d({\bf g}_{n+1})= r({\bf  f}_{n+1}) = r({\bf f}_n) = r({\bf g}_{n+1})$ and $\infty>e(Y_{n+1})=e(Y_n)$;  an induction argument together with Theorem~\ref{thm:lifting}\eqref{p4} shows that ${\bf g}_{n+1}$  preserves $0$--stabilizers for $n\geq 1$.

The third statement is proved as follows. First observe that
\[\infty>d({\bf f})= d({\bf f}_{n+1})\geq r({\bf f}_{n+1})= d({\bf g}_{n+1}) \geq r({\bf g}_{n+1})\geq r({\bf f}_n)\geq 0,\]
and hence 
\[ d({\bf f})-r({\bf f}_n) \geq d({\bf f})-r({\bf f}_{n+1})\geq 0.\]
If $r({\bf f}_n) = r({\bf  f}_{n+1})$ then $d({\bf g}_{n+1})= r({\bf g}_{n+1})$. Then Lemma~\ref{lem:isomorphism}\eqref{lem:t2r2} implies that ${\bf g}_{n+1}\colon Y_{n+1}\to Y_n$ is a $G$--equivariant inclusion and therefore $e(Y_{n}) \geq e(Y_{n+1})$. It follows that $|{\bf f}_{n+1}| \leq_{dic} |{\bf f}_{n}|$.\end{proof}

To conclude the proof of  Theorem~\ref{thm:maximal-towers-full} observe that if for every $n$ the complex $Y_n$ is not simply-connected, then there is an infinite strictly decreasing sequence of elements of $(\N\times \N, \leq_{dic})$. Since this is impossible, there is $n$ such that $Y_n$ is simply-connected and then ${\bf f}_n$ is a maximal $\mathcal F$--tower lifting of $\bf f$. 
\end{proof}

\subsection{Existence of Maximal Towers Liftings}

\begin{theorem}[Maximal Equivariant Towers]\label{thm:maximal-towers}
Let ${\bf f}\colon X \to Y$ be a $G$--map. If $X$ is one-connected and $G$--cocompact, then $\bf f$ has a maximal equivariant tower lifting.
\end{theorem}

\begin{proof}[Sketch of the proof]
The proof of Theorem~\ref{thm:maximal-towers} is the same as the proof of Theorem~\ref{thm:maximal-towers-full} with some simplifications.
Define inductively  the sequence of tower liftings ${\bf f}_{n+1} \colon X \to Y_{n+1}$ where $Y_n$ is defined as the $G$--subcomplex ${\bf \tilde f}_{n}(X)$ of $\widetilde Y_n$ and $Y_1=f(X)$. Since $X$ is $G$--cocompact and ${f}_n$ is surjective, $Y_n$ is $G$--cocompact and hence $d({\bf f}_n) - r({\bf f}_n)<\infty$. Then, the same argument, shows that  
\[d({\bf f}_n)-r({\bf f}_n) \geq d({\bf f}_{n+1})-r({\bf f}_{n+1}).\] One verifies that if $Y_n$ is not simply-connected then the inequality above is strict. Specifically, if 
$r({\bf f}_n) = r({\bf f}_{n+1})$ then $g_{n+1}$ is a surjective $G$--equivariant tower preserving $0$--stabilizers and such that $d({\bf g}_{n+1})=r({\bf g}_{n+1})$; then  Lemma~\ref{lem:isomorphism}\eqref{lem:third} implies that $g_{n+1}$ is an isomorphism factoring through $\widetilde Y_n$ and hence $Y_n$ is simply-connected. If each $Y_n$ is not simply-connected then $d({\bf f}_n)-r({\bf f}_n)$ defines a strictly decreasing infinite sequence of natural numbers which is impossible. Therefore some $Y_n$ is simply-connected and the corresponding ${\bf f}_n$ is the desired maximal tower lifting of ${\bf f}$.
\end{proof}

\section{Immersed Cocompact Cores}\label{sec:cores}

\begin{theorem}\label{thm:main-full}  If $Y$ is a one-connected, proper and locally finite  $H$--complex and $G\leq H$ is finitely presented, then there is a one-connected cocompact $G$--complex $X$ and a  $G$--equivariant $\mathcal F$--tower $X\to Y$.
\end{theorem}

The proof of the theorem requires two lemmas.  

\begin{lemma}\label{lem:last-one}
Let $\mathcal G$ be a finite graph of groups such that vertex groups are finitely generated and edge groups are finite. If $\varphi\colon \pi_1 (\mathcal G) \to G$ is a surjective group homomorphism into a finitely presented group, then $\kernel(\varphi)$ is normally finitely generated.
\end{lemma}
\begin{proof}
The hypotheses on $\mathcal G$ imply that   $\pi_1 (\mathcal G)$ is finitely generated. Therefore there is a surjective homomorphism $F \overset{\psi}{\to} \pi_1 (\mathcal G)$  where $F$ is a finite rank free group. Since $G$ is finitely presented and $F$ has finite rank, the kernel of $\varphi\circ \psi$ is normally finitely generated, say $\kernel (\varphi \circ \psi)=\nclose{r_1, \ldots , r_m}$.  Since $\kernel ( \varphi ) = \psi \left ( \kernel (\varphi \circ \psi) \right)$, we have that $\kernel (\varphi)=\nclose{\psi(r_1), \ldots , \psi(r_m)}$.
\end{proof}

Recall that a group is \emph{slender} if all its subgroups are finitely generated, and a $G$--complex is \emph{almost proper} if $G$ acts properly on the complement of the $0$--skeleton.

\begin{lemma}[One-connected Complex for Finitely Presented Subgroup] \label{lem:map} Let $Y$ be a one-connected, almost proper   $H$--complex such that $H$--stabilizers of cells are slender.  Suppose that $G\leq H$ is finitely presented. Then there exists a cocompact and one-connected \Gtwo--complex $X$, and an equivariant map ${\bf f} \colon X \rightarrow Y$ such that $f_\#$ is the inclusion $G \hookrightarrow H$ and $f_\#$ is injective when restricted to $0$--cell stabilizers. 
\end{lemma}
\begin{proof}
First we construct the $1$--skeleton $X^{(1)}$ of $X$ as a $G$--equivariant cocompact subcomplex of $Y$. Let $\{g_i\}_{i=1}^m$ be a finite generating set for $G$ and let  $y_0$ be a $0$--cell of $Y$. Since $Y$ is connected, for each $1\leq i \leq m$ there is a combinatorial path $\gamma_i$ from $y_0$ to $g_iy_0$. Let $D$ be the finite subcomplex $D=\gamma_1 \cup \cdots \cup \gamma_m$ of $Y$ and let $X^{(1)}$ be the union of all $G$--translates of $D$ in $Y$.  Then $X^{(1)}$ is a  $G$--equivariant $1$-dimensional cocompact connected  subcomplex of $Y$.  Connectedness follows from the assumption that  $\{g_i\}_{i=1}^m$ generates $G$, and cocompactness from $D$ being a finite subcomplex.  

Next we show that $\pi_1   X^{(1)} $ is a normally finitely generated group. Invoking Theorem~\ref{thm:lifting},  consider the group $\widetilde G$ acting on the universal cover  $\mathcal T$ of $X^{(1)}$ such that $\mathcal T  \to X^{(1)}$ is equivariant with respect to $\widetilde G \to G$  and  $\pi_1 X^{(1)}$  is isomorphic to the kernel of $\widetilde G \to G$.    By considering the barycentric subdivision $\mathcal T'$ of $\mathcal T$, we have a cocompact action without inversions of $\widetilde G$ on the tree $\mathcal T'$; here cocompactness follows from $X^{(1)}$ being $G$--cocompact.  By Theorem~\ref{thm:lifting}~\eqref{p4}  the group homomorphism  $g_\#\colon \widetilde G_x \to G_{g(x)}$ is injective for every $x\in \mathcal T$. Therefore  $1$--cell $\widetilde G$--stabilizers of $\mathcal T'$ are finite, since $Y$ is an almost proper $G$--complex. Moreover,  $0$--cell $\widetilde G$--stabilizers of $\mathcal T'$ are  isomorphic to either a subgroup of the $G$--stabilizer of a $0$--cell of $Y$ and hence finitely generated by the slender hypothesis; or to a subgroup of the setwise $\widetilde G$--stabilizer of a $1$--cell of $Y$ and hence finite since $Y$ is almost proper.  Then the theory of Bass and Serre on actions on trees~\cite{Serre-Trees} implies that  $\widetilde G$ is isomorphic to the fundamental group of a finite graph of groups $\mathcal G$ with finite edge groups and finitely generated vertex groups. By Lemma~\ref{lem:last-one}, the kernel of $\widetilde G \to G$ is normally finitely generated. 

To conclude the proof, we paste finitely many $G$-orbits of $2$-cells to $X^{(1)}$ to obtain a one-connected complex. Choose a $0$--cell $\tilde x_0$ of $\mathcal T$ as a basepoint. Since $\pi_1 X^{(1)}$ is normally finitely generated, there is a finite collection $\{ r_i\}_{i=1}^q$ of based loops in $X^{(1)}$ such that  $\pi_1 X^{(1)}= \langle \langle r_1, \dots, r_q \rangle \rangle$.  Since $Y$ is simply-connected and $X^{(1)}$ is a subcomplex of $Y$, for each $r_i$ there is   a disk-diagram $D_i \to Y$ with boundary path $r_i$.
Let $X$ be the complex obtained by attaching a copy of $D_i$ to $X^{(1)}$ along the closed path $g.r_i$  for each $g\in G$ and each $1\leq i \leq q$. Observe that $X$ is connected and simply-connected, and the $G$-action on $X^{(1)}$ naturally extends to a cocompact  $G$-action on $X$. The equivariant inclusion $X^{(1)} \hookrightarrow Y$ extends to a map $X \to Y$ equivariant with respect to $G\hookrightarrow H$. 
\end{proof}

The proof of the Theorem~\ref{thm:main-full} is an application of  Lemma~\ref{lem:map} together with Theorem~\ref{thm:maximal-towers-full}.
\begin{proof}[Proof of Theorem~\ref{thm:main-full}]
Observe that the hypotheses of Lemma~\ref{lem:map} are satisfied since $Y$ being proper implies that stabilizers of $0$--cells are finite and proper implies almost proper. Therefore there is a one-connected cocompact $G$--complex $X_0$ and a  $G$--map ${\bf f}\colon X_0 \to Y$. By Theorem~\ref{thm:maximal-towers-full}, there is a maximal $\mathcal F$--tower lifting $\bf {f}= \bf{g} \circ \bf{f'}$ of $\bf f$, where $\mathbf{g}\colon X\to Y$. By Proposition~\ref{prop:maxitower}~\eqref{item:maxitower-full}, $f'$ is $0$--surjective and $\pi_1$--surjective, and $f'_\#$ is surjective. It follows that $X$ is one-connected and $f'_\#$ is an isomorphism, in particular ${\bf f'}$ is a $G$--map. Since  $Y$ is locally finite and ${\bf g}$ is a tower, it follows that $X$ is locally finite.  Since $X_0$ is $G$--cocompact and $f'\colon X_0\to X$ is $0$--surjective, we have that $X$ is $G$--cocompact.
\end{proof}
Analogously to the argument in the proof of Theorem~\ref{thm:main-full}, one obtains the following result by combining Lemma~\ref{lem:map} and Theorem~\ref{thm:maximal-towers}. In this result,  the complex  $Y$ is not necessarily locally finite but the conclusion is weaker. 
\begin{theorem}\label{thm:main} 
Let $Y$ be a one-connected and almost proper $H$--complex such that $H$--stabilizers of cells are slender. If $G\leq H$ is finitely presented, then there is a one-connected cocompact  $G$--complex $X$ and a $G$--equivariant tower $X\to Y$.
\end{theorem}
A version of Theorem~\ref{thm:main} appears in~\cite[Lem. 6.4]{MaWi13} where is shown that a $G$--map with locally finite target  always factors as the composition of a surjective and $\pi_1$--surjective $G$--map followed by a $G$--equivariant immersion.

\section{Applications}\label{sec:applications}

\subsection{Subgroups of Diagramatically reducible groups}

\begin{definition}[Near-immersion]
A map $X\to Y$ is a \emph{near immersion} if it is locally injective in the complement of the $0$--skeleton of $X$. 
\end{definition}

\begin{definition}[Diagramatically Reducible Complex]\cite{Ge87}
A $2$--complex $X$ is \emph{diagramatically  reducible} if there are no near-immersions $C\to X$, where  $C$ is a cell structure for the $2$--dimensional sphere.
\end{definition}

First we recall some properties of diagramatically reducible complexes in the proposition below.

\begin{proposition}\label{thm:DRprop}
\begin{enumerate}
\item[ ]
\item \label{DRaspherical} Diagramatically reducible complexes are aspherical.
\item \label{DRtowers} Covers and subcomplexes of diagramatically reducible complexes are diagramatically reducible.
\item \label{DRsubdivision} The barycentric subdivision of a diagramatically reducible complex is diagramatically reducible. 
\end{enumerate}
\end{proposition}
\begin{proof}
The first statement is a result of Gersten~\cite[Rem. 3.2]{Ge87}, the second statement is trivial since the composition of a near-immersion and an immersion is a near-immersion, and the third statement is due to Howie~\cite[Rem. 6.10]{Ge87}.
\end{proof}

\begin{definition}[Free $1$--cells and Collapsing]
Let $X$ be a $2$--complex. A $1$--cell $e$ of a subcomplex $Z \leq X$ is \emph{free} if it belongs to the boundary of a $2$--cell $f$ of $Z$, and $e$ does not belong to the boundary of a $2$--cell $\neq f$ in $Z$. In this case, \emph{collapsing $Z$ along $e$} means to remove the interior of $e$ and the interior of $f$. 
\end{definition}

The following characterization of diagramatically reducible complex is a result of  Jon Corson. 
\begin{theorem}[Characterization]\cite[Thm. 2.1]{Co01}\label{thm:DRchar}
A one-connected $2$--complex is diagramatically reducible if and only if every finite subcomplex is $1$--dimensional or contains a free $1$--cell.
\end{theorem}

\begin{remark}[Equivariant Collapsing and Inversions]\label{rem:DRequiv-collap}
Recall that a group action on a complex has no inversions if whenever a cell is fixed setwise by a group element then it is fixed pointwise by the group element.  Let $Z$ be a $G$--complex without inversions and suppose that $e$ is a free $1$--cell of $Z$ that belongs to the boundary of the $2$--cell $f$. Observe that for every $g\in G$ the $1$--cell $g.e$ is free in $Z$.  Since $G$ acts without inversions, for every $g\in G$, the $2$--cell $g.f$ contains only one $1$--cell in the $G$--orbit of $e$, namely, $g.e$. Therefore we can simultaneously collapse $Z$ along $g.e$ for every $g\in G$ obtaining a $G$--equivariant subcomplex $Z'$ of $X$. 
\end{remark}

Corson also proves that if $F$ is a finite group acting on a one-connected diagramatically reducible $2$--complex $X$ then the fixed point set $X^F$ of $F$ is non-empty~\cite[Thm. 4.1]{Co01}. The following proposition shows that $X^F$ is also contractible. 
\begin{proposition}[Contractible Fixed Point Sets]\label{prop:contractible}
Let $X$ be a one-connected diagramatically reducible $2$--complex. If $F$ is a finite group acting on $X$ without inversions, then the fixed point set $X^F$ of $F$ is a non-empty contractible subcomplex.
\end{proposition}
\begin{proof}
Since $F$ acts without inversions, $X^F$ is a subcomplex of $X$.  By Proposition~\ref{thm:DRprop}\eqref{DRaspherical}, it is enough to show that $X^F$ is one-connected.  First we verify that $X^F$ is connected.  Take two vertices $x_0$ and $x_1$ of $X^F$. Since $X$ is connected,  there is an edge--path  $\alpha$  in $X$ between $x_0$ and $x_1$.  Let $Y_0$ be the subcomplex of $X$ defined as  the union of all the images of $\alpha$ under the action of $F$. Observe that $Y_0$ is a connected finite $1$-dimensional subcomplex of $X$ invariant under the $F$-action.  Construct a finite one-connected $F$-complex $Y$ as follows. Fix a basepoint of $Y_0$ and let $\gamma_1, \ldots , \gamma_n$ be a collection of closed paths in $Y_0$ that generate $\pi_1 (Y_0)$. For each $\gamma_i$, there is a disk diagram $D_i \to X$ with boundary path $\gamma_i \to X$. Let $Y$ be the $2$--complex obtained by attaching a copy of $D_i$ to $Y_0$ along $k.\gamma_i$ for each $k\in F$. Then $Y$ is a finite one-connected $2$--complex, the $F$-action on $Y_0$ extends to an action on $Y$, and there is a natural $F$-equivariant map $f \colon Y\to X$.   By Theorem~\ref{thm:maximal-towers}, there is a maximal equivariant tower lifting $Y\overset{f'}\to Z \overset{g}\to X$ of $f$.  
Since $g$ is an immersion, $Z$ is diagramatically reducible.  Since $F$ acts without inversions on $X$, it also acts without inversions on $Z$.  By Theorem~\ref{thm:DRchar}, it follows that if $Z$ contains $2$--cells then it has a free $1$--cell $e$.  After a finite number of $F$-equivariant collapses of $Z$ one obtains a $1$-dimensional one-connected $F$-complex equivariantly immersed into $X$, see Remark~\ref{rem:DRequiv-collap}.  Without loss of generality, we can assume that $Z$ is $1$-dimensional. Then $Z$ is a tree and therefore $F$ fixes pointwise an edge-path in $Z$ between $x_0$ and $x_1$. By equivariance of the map $Y\to X$, $F$ fixes pointwise a path between $x_0$ and $x_1$.

Now we verify that $X^F$ is simply-connected. Since $X$ is simply-connected, consider an essential embedded closed path $\gamma$ in $X^F$ with minimal area in $X$. Then there is a disk-diagram $D\to X$ with boundary $\gamma$ of minimal area. Let $g\in F$ with $g\neq 1$ and observe that the spherical diagram $D\cup gD \to X$ is a near-immersion. Since this is impossible, $X^F$ is simply connected.
\end{proof}

\begin{remark}[Inversions and Connected Fixed Point Sets]
During the review process of the article, the referee observed that in the proof of Proposition~\ref{prop:contractible}, one can prove that $X^F$ is connected without assuming that $F$ acts without inversions.  The sketch of the argument is as follows.  Suppose that $X^F$ is not connected and choose  $x_0$ and $x_1$ at minimal distance in different connected components of $X^F$. Consider a path $\gamma$ in $X$ between $x_0$ and $x_1$. Since $x_0$ and $x_1$ are in different connected components, the path $\gamma$ is not fixed by $F$ and hence there is a non-trivial element $g\in F$ that does not fix $\gamma$ pointwise.   For given $\gamma$ and $g\in F$,  there is a disk diagram $D \to X$ with boundary path $\gamma^{-1} g(\gamma)$. Among all these possible choices of $\gamma$, $g$ and $D$,  choose the ones that minimize $\area (D)$.   It follows that the diagram $D$ contains $2$--cells and no cut-points.  Let $n$ be the order of the element $g$. Then one can glue together  $n$ copies of $D$, by identifying $g(\gamma)$ in the $i$-th copy of $D$ with $\gamma$ in the $i+1$-copy of $D$, producing a sphere $S$. The map $S\to X$ which maps the $i$-th copy of $D$ in $S$ to $g^i(D)$ in $X$ is a near-immersion by our minimality choices. This contradicts that $X$ is diagramatically reducible. 
\end{remark}

\begin{definition}
Let $G$ be a group. A proper $G$--complex $X$ is a model for $\underline{E}G$ if for every finite subgroup of $F\leq G$, the fixed point set $X^F$ is contractible.
\end{definition}

\begin{theorem}[Diagramatically Reduced Groups]\label{thm:EG}
Let $Y$ be a diagramatically reducible proper \Htwo--complex. If $G\leq H$ is finitely presented then $G$ admits a diagramatically reduced  $2$--dimensional  cocompact model for $\underline{E}H$.
\end{theorem}
\begin{proof}
By passing to a subdivision of $Y$ we can assume that $G$ acts without inversions while still assuming that $Y$ is diagramatically reducible, see Proposition~\ref{thm:DRprop}\eqref{DRsubdivision}.   By Theorem~\ref{thm:main}, there is an equivariant immersion $X\to Y$ where $X$ is a one-connected cocompact \Gtwo--complex. Since $Y$ is diagramatically reducible, Proposition~\ref{thm:DRprop}\eqref{DRtowers} implies that $X$ is diagramatically reducible as well. Since the $H$-action on $Y$ is proper and without inversions, the same properties hold for the $G$--action on $X$.   Let $K$ be a finite subgroup of $G$. Then Proposition~\ref{prop:contractible} implies that the fixed point set $X^F$ is contractible. 
\end{proof}

\subsection{Subgroups of Relatively Hyperbolic Groups Acting on Negatively Curved $2$-Complexes}\label{subsec:relhyp}  

\begin{definition}[Fine Graphs and Fine Complexes~\cite{Bo12}]
A $1$-complex is \emph{fine} if each $1$--cell is contained in only finitely many circuits of length $n$ for each $n$. Equivalently,  the number of embedded paths of length $n$ between any pair of (distinct) $0$--cells is finite.  A complex is \emph{fine} if its $1$--skeleton is fine.
\end{definition}

\begin{definition}[Relatively Hyperbolic Groups~\cite{Bo12}]\label{def:rh}
A group $G$ is \emph{hyperbolic relative to a finite collection of subgroups $\mathcal P$} if $G$ acts cocompactly, almost properly on a connected, fine, $\delta$--hyperbolic $1$--complex, and $\mathcal P$ is a set of representatives of distinct conjugacy classes of vertex stabilizers  such that each infinite stabilizer is represented. 
\end{definition}

\begin{proposition}[$2$--dimensional Relative Hyperbolicity]\label{prop:relhyp}
Let $X$ be a one-connected, negatively curved, fine, cocompact, and almost proper \Gtwo--complex. Then $G$ is a hyperbolic group relative to a (hence any) collection of representatives of conjugacy classes of $0$--cell stabilizers.
\end{proposition}
\begin{proof}
Since $X$ is negatively curved and one-connected, it is a $\delta$-hyperbolic space. Indeed, it is well known that a $CAT(\kappa)$--space with $\kappa<0$ is $\delta$--hyperbolic, and in the conformal case  $X$ satisfies a linear isoperimetric inequality and hence the combinatorial metric on its one skeleton is a $\delta$--hyperbolic space~\cite{BrHa99, Ge96}. Since $X$ is $G$--cocompact there are finitely many types of $2$--cells and hence the $1$--skeleton $X^{(1)}$ is quasi-isometric to $X$; in particular $X^{(1)}$ is a $\delta'$-hyperbolic space.  It follows that $X^{(1)}$ is endowed with a $G$--action satisfying the conditions of Definition~\ref{def:rh}.
\end{proof}

\begin{figure}
\labellist
\small\hair 2pt
\pinlabel $x_0$ at 180 180
\pinlabel $x_1$ at 518 166
\pinlabel $\eta$  at 292 356
\pinlabel $a$ at 15 235
\pinlabel $b$ at 215 5
\pinlabel $e_1$ at 99 200
\pinlabel $e_2$ at 191 101
\pinlabel $D$ at 609 376
\endlabellist
\centering
\includegraphics[width=0.7\linewidth]{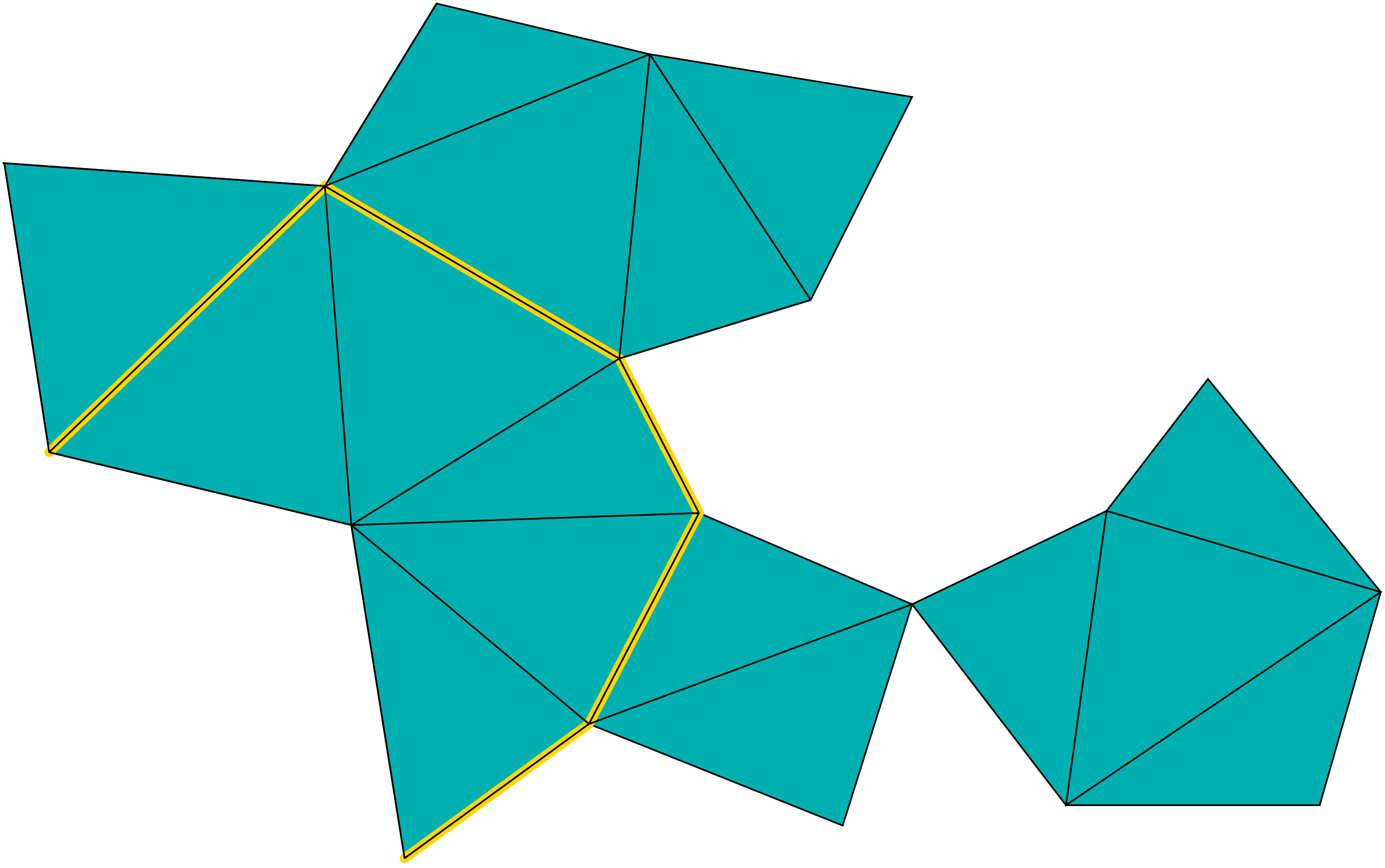}
\caption{In the disk diagram $D$, the space $D\setminus \{x_1\}$ is not connected.  If $D$ is a disk diagram, $x_0$ is a $0$--cell on $\partial D$ adjacent to the $0$--cells $a,b \in \partial D$ with $a\neq b$, and  $D\setminus\{x_0\}$ is connected, then there is an edge  path $\eta \to D$ between $a,b$ that factors through $D\setminus \{x_0\}$ and goes around the $2$--cells  of $D$ adjacent to $x_0$. The length of $\eta$ is bounded by $\area (D)\cdot C_1$ where $C_1$ is an upper bound for the boundary length of $2$--cells  of $D$. 
 }
\label{fig:gallery}
\end{figure}

\begin{proposition}\label{prop:fine-is-star-2} \label{prop:fine-is-star}
Let $Y$ be a diagramatically reducible one-connected cocompact almost proper $\Htwo$-complex with fine one-skeleton. If $X$ is one-connected and there is an immersion $X\to Y$ then $X$ has fine one-skeleton.
\end{proposition}
\begin{proof}
We use the following characterization of simplicial fine graphs due to Brian Bowditch~\cite[Prop. 2.1]{Bo12}. By simplicial we mean no double edges and no single edge loops.  A simplicial graph $K$ is fine if and only if for each vertex $x\in K$, the set $V(x)$ of vertices adjacent to $x$ has the following property: every subset of $V(x)$ which is  bounded in $K\setminus \{x\}$ with respect to the combinatorial metric is finite.  

Without loss of generality assume that the boundary path of every $2$--cell of $Y$ is an embedded path, and that $X$ and $Y$ have simplicial one-skeleton.   Indeed, by considering the barycentric subdivisions of $X$ and $Y$, we can assume the boundary paths of $2$--cells are embedded and one-skeletons are simplicial. Proving the proposition for the barycentric subdivisions is sufficient since the one-skeleton of an almost proper $\Htwo$-complex is fine if and only if the one-skeleton of its  barycentric subdivision is fine, this follows directly from~\cite[Lem. 2.3]{Bo12} or~\cite[Lem. 2.9]{MaWi11b}; moreover,  a complex is diagramatically reducible if and only if its barycentric subdivision is diagramatically reducible, see Proposition~\ref{thm:DRprop}\eqref{DRsubdivision}. 

Since $Y$ admits a cocompact $H$--action and its one skeleton is fine, for every $n\in \N$ there are finitely many cycles of length $n$ up to the $H$--action.
Therefore, since $Y$ is simply-connected, there is a well defined Dehn function $\Delta \colon \N \to \N$, that is, $\Delta (n)$ is an upper bound for the area of minimal area disk diagrams with  given boundary path  of length $\leq n$.

Let $f\colon X\to Y$ be an immersion. Let $x_0\in X$  be $0$--cell of $X$ and let $A$ be a subset of $0$--cells adjacent to $x$. Denote by  $\diam (A)$ the diameter of $A$ in $X\setminus \{x_0\}$ in the combinatorial metric, and analogously let $\diam (f(A))$ denote the diameter of $f(A)$ in $Y\setminus \{f(x_0)\}$ in the combinatorial metric. The claim is that if $\diam (A)$ is finite then
\begin{equation} \label{eq:fine} \diam \left (f(A) \right) \leq C \cdot  \Delta \left( \diam (A) +2 \right),\end{equation} where $C$ is an upper bound for the boundary length of a $2$--cell of $Y$; here $C$ is finite since $Y$ admits a cocompact action.    Assuming the claim, we conclude using Bowditch's characterization of fineness as follows.  If $\diam (A)$ is finite, then the claim implies that  $\diam (f(A))$ is finite;  then $Y$ being fine implies that  $f(A)$ is a finite set; since $f$ is an immersion and one-skeletons are simplicial, the induced map $f\colon A \to f(A)$ is a bijection and hence $A$ is finite set. 

Suppose that $\diam (A)=m<\infty$ in $X\setminus \{x_0\}$. Let $a,b\in A$ and suppose that $a\neq b$. Then there is a combinatorial path $\gamma$ in $X\setminus \{x_0\}$ from $a$ to $b$ of length $\leq m$. If the path $f\circ \gamma$ in $Y$ does not contain the $0$--cell $f(x_0)$ then the combinatorial distance between $f(a)$ and $f(b)$ is bounded by $m$; however this assumption on $f\circ \gamma$ might not hold.   A general argument is as follows.  Consider the closed path $\gamma'\to X$ defined as the concatenation  $\gamma'=e_1 \ast \gamma \ast e_2$ where $e_1$ is an $1$--cell from $x_0$ to $a$, and $e_2$ is a $1$--cell  from $b$ to $x_0$.  Since $X$ is simply-connected, there is a near-immersion of a disk-diagram $D\to X$ with boundary path $\gamma'$. Observe that there is only one $0$--cell in $\partial D$ mapping to $x_0 \in X$; by abuse of notation,  let $x_0$ denote this $0$--cell of $D$. 

The main observation is that $D\setminus \{x_0\}$ is connected. Indeed, if $D\setminus \{x_0\}$ is not connected, then the boundary path $\partial D \to X$ pass through $x_0$ more than once.  Since  $e_1\to X$ and $e_2\to X$ are $1$--cells with only one endpoint equal $x_0$, and the image of $\gamma \to X$ does not contain $x_0$, it follows that $\partial D = e_1 \ast \gamma \ast e_2$ pass through $x_0$ only once, and hence $D\setminus \{x_0\}$ is connected.

The fact that $D\setminus \{x_0\}$ is  connected implies that there is an embedded path $\eta \to D$ between  the two $0$--cells adjacent to $x_0$ in $\partial D$,  this path $\eta\to D$ factors through $D\setminus \{x_0\}$ and goes around the $2$--cells  of $D$ adjacent to $x_0$; see Figure~\ref{fig:gallery}.  Observe that $\eta \to D \to X \overset{f}\to Y$ is a path between $f(a)$ and $f(b)$, and  the combinatorial length of $\eta$ is bounded by $C \cdot \area (D)$ where $C$ is the upper bound for the boundary length of a $2$--cell of $Y$.

Now observe that  the path $\eta \to D \to X \overset{f}{\to} Y$ does not intersect $f(x_0)$. Indeed, since boundary paths of $2$--cells  of $Y$ (and hence of $D$) are embedded, if $\eta \to Y$ intersects $f(x_0)$ then there is a $2$--cell $R$ of $D$ whose boundary path contains $x_0\in \partial D$ and another $0$--cell $x \in \eta$ both mapping to $f(x_0)$; this would imply that $\partial R \to Y$ is not an embedded path which is impossible by our initial assumption. 

Since $Y$ is diagramatically reducible and  $D\to Y$ is a near-immersion, it follows that $D\to Y$ is a minimal area disk diagram for $\partial D \to Y$. Therefore $|\eta| \leq C \cdot \Delta (|\partial D|)$, and hence the combinatorial distance between $f(a)$ and $f(b)$ in $Y\setminus \{ f(x_0)\}$ is bounded by $|\eta| \leq C  \cdot  \Delta \left( \diam (A) +2 \right)$. Since $a$ and $b$ were arbitrary, we have proved that inequality~\eqref{eq:fine} holds.
\end{proof}


\begin{theorem}[Subgroups of $2$-Dim. Rel. Hyp. Groups are Rel. Hyp]\label{thm:rel-hyp}
Let $Y$ be a one-connected negatively curved, fine,  almost proper and cocompact \Htwo--complex such that  $H$--stabilizers of cells are slender.  If $G\leq H$ is finitely presented then $G$ is  hyperbolic relative to a finite collection of $G$--stabilizers of cells of $Y$.
\end{theorem}
\begin{proof}[Proof of Theorem~\ref{thm:rel-hyp}]
By Theorem~\ref{thm:main}, there is an equivariant tower $X\to Y$ where $X$ is a one-connected cocompact \Gtwo--complex.  Since $Y$ is negatively curved and almost proper $G$--complex, the same properties hold for $X$.   Proposition~\ref{prop:fine-is-star} implies that $X$ is fine. 
Then Proposition~\ref{prop:relhyp} implies that $G$ is hyperbolic relative to a finite collection $\mathcal P$ of $G$--stabilizers of cells of $X$. Since $X \to Y$ is a $G$--map, the collection $\mathcal P$ consists of $G$--stabilizers of cells of $Y$. 
\end{proof}

%
%
%
\bibliographystyle{gtart}

\end{document}